\documentclass[a4paper,10pt]{article}

\usepackage{relsize,exscale}
\usepackage{color}

\usepackage{blkarray}
\usepackage[arrow,matrix]{xy}
\usepackage{hyperref}
\usepackage{makeidx}
\usepackage{bm}
\usepackage{latexsym}
\usepackage{amscd}
\usepackage{amsmath}
\usepackage{amssymb}
\usepackage{hhline}
\usepackage{mathenv}
\usepackage{stmaryrd}

\usepackage{amsthm}
\usepackage{float}
\usepackage{graphicx}
\usepackage{psfrag}

\theoremstyle{plain}
\newtheorem{theorem}{Theorem}[section]

\theoremstyle{definition}

\newtheorem{remark}[theorem]{Remark}

\newcommand{\Q}{\mathbb{Q}}

\newcommand{\C}{\mathbb{C}}

\newcommand{\PP}{\mathbb{P}}

\newcommand{\X}{{\mathcal X}}

\newcommand{\CC}{{\mathcal C}}
\newcommand{\A}{{\mathcal A}}

\newcommand{\PPP}{{\mathcal P}}

\definecolor{deepblue}{cmyk}{0,0.83,1,0.70}
\definecolor{gray}{cmyk}{0,0,0,0.3}
\definecolor{rred}{cmyk}{0,1,1,0}
\definecolor{chairo}{cmyk}{0,0.83,1,0.70}
\definecolor{roypur}{cmyk}{0.75,0.90,0,0.1}
\definecolor{darkorc}{cmyk}{0.40,0.80,0.20,0}
\definecolor{oliv}{cmyk}{0.64,0.00,0.75,0.56}
\definecolor{azuro}{cmyk}{1,1,0,0.46}

\usepackage{color}
\usepackage{ifthen}

\newboolean{usecolor}
\setboolean{usecolor}{true}

\ifthenelse{\boolean{usecolor}}
{
  \definecolor{colore}{cmyk}{0,1,0.6,0}
  \definecolor{coloregen}{cmyk}{0.7,0,1,0}
  \definecolor{colorealex}{cmyk}{1,0.6,0,0}
}
{
  \definecolor{colore}{cmyk}{0,0,0,1}
  \definecolor{coloregen}{cmyk}{0,0,0,1}
  \definecolor{colorealex}{cmyk}{0,0,0,1}
}

\title{ON TWO PENCILS OF CUBIC CURVES
}
\author{Pauline Bailet\thanks{Department of Mathematics and Computer Science, Universit\"{a}t Bremen, 28359 Bremen, Germany. Email: pauline.bailet@uni-bremen.de}}
\date{\today}
\begin{document}
\maketitle

\begin{abstract}
We give two examples of curve arrangements $\CC: f=0 \subset \PP^2_\C$ of pencil type which are very close to line arrangements, though the action of the monodromy $h^1$ on the cohomology $H^1(F,\C)$ of the Milnor Fiber $F:= f^{-1}(1) \subset \C^3$ has eigenvalues of order 5 and 6, showing that surprising situations can occur for larger classes of curve arrangements than for line arrangements. Our computations rely on the algorithm given by A. Dimca and G. Sticlaru in \cite{D3} which detects the non trivial monodromy eigenspaces of free curves. 
\end{abstract}
\noindent
{\bf MSC-class}: 14H50, 14DXX, 14FXX\\
{\bf Keywords}: Pencil of curves, Curve arrangements, Milnor fiber, Monodromy.

\section{Introduction}

Let $\mathcal{C} : f = 0$ be a reduced curve of degree $d$ in the complex projective plane $\PP^2_\C.$ Consider the corresponding complement $U = \PP^2_\C \backslash \CC,$ and the global Milnor fiber $F$ defined by $f(x,y,z) = 1$ in $\C^3$ with monodromy action $h: F \rightarrow F,\, h(x) = \exp(2\pi i/d)\cdot(x,y,z).$ To determine the eigenvalues of the monodromy operator 
\begin{equation}\label{eq1}
h^1 : H^1(F,\C) \rightarrow H^1(F,\C)
\end{equation}  is a rather difficult problem. Assume that $\CC$ is a curve arrangement of pencil type, i.e. the defining equation has the form $$f = q_1q_2 \cdots q_m,$$ for some $m \geq 3,$ where $\deg q_1 = \cdots = \deg q_m = k \geq 2$ and the curves $\CC_i : q_i = 0$ for $i = 1,\hdots,m$ are all members of the pencil $\PPP : u\CC_1 + v\CC_2.$ The situation when the curves $\CC_i$ are line arrangements was systematically considered, see \cite{Falk}, \cite{PS}, \cite{Yuz}, particularly with a view in understanding whether the monodromy action (\ref{eq1}) is combinatorially determined. We mention just three striking facts in this direction, assuming that $\CC_i$ are line arrangements. 
\begin{enumerate}
\item The number $m$ of members of the pencil $\PPP$ is at most 4, and the Hessian arrangement is the only known such pencil type arrangement with $m = 4,$ see \cite{Yuz}. For this arrangement the corresponding eigenvalues are the roots of unity of order 4. 
\item If the line arrangement $\CC$ has only double and triple points, then the monodromy operator (\ref{eq1}) is combinatorially determined, and the corresponding eigenvalues are cubic roots of unity, see \cite{PS}.
\item It is not known whether the monodromy operator (\ref{eq1}) can have eigenvalues which are not roots of unity of order 3 or 4, see \cite{PS}. 
\end{enumerate} 
On the other hand, it is known that if we consider larger classes of curve arrangements, e.g. conic and line arrangements, many new properties can occur. For instance, Terao's conjecture about the combinatorial invariance of freeness is open for line arrangements, but false for conic and line arrangements, see \cite{ST}. Our main result, to be stated next, says that a similar situation occurs when looking at the properties 1. and 3. listed above. 

\begin{theorem}\label{thm1}
Consider the pencil of cubic curves given by $\PPP : u\CC_1 + v\CC_2$ in the complex projective plane $\PP^2_\C,$ where 
\begin{center}
$\CC_1 : q_1(x,y,z) = 3xyz + y^3 + z^3 = 0$ and $\CC_2 : q_2(x,y,z) = 3xyz + x^3 + z^3 = 0.$
\end{center}

Then the following hold. 
\begin{enumerate}
\item The two curves $\CC_1$ and $\CC_2$ meet transversally in 9 points, hence the generic member of the pencil $\PPP$ is smooth.
\item The number $m$ of singular members in the pencil $\PPP$ is 5 and they are the fibers of the points $(u':v')\in \PP_\C^1$ that are listed as follows: for $(u' : v') = (1 : 0)$ and $(u' : v') = (0 : 1)$ we get two nodal irreducible cubics, for $(u' : v') = (1 : 1)$ we get a union of three concurrent lines, and for $(u' : v') = (1:-t_i),$ where $i = 1,2$ and $t_i$'s are the roots of the equation $$t^2 + 3t + 1 = 0,$$ we get each time a triangle.
\item The union of the singular members in the pencil $\PPP$ is a free curve $$\CC : f = q_1q_2(q_1- q_2)(q_1^2-3q_1q_2 + q_2^2) = 0$$ of degree $d=15$ and with exponents $d_1 = 4$ and $d_2 = 10,$ having 9 lines and two nodal cubics as irreducible components.
\item The monodromy operator (\ref{eq1}) associated to the curve $\CC$ has eigenvalues which are roots of unity of order 5. More precisely, the Alexander polynomial of the curve $\CC$ is given by $$\Delta_\CC (t) = (t-1)^7(t^5 -1)^3.$$ 
\end{enumerate}
\end{theorem}
  
\begin{remark} \label{rk1}
The curve arrangement $\CC$ in Theorem \ref{thm1} contains the line arrangement $$\A : (q_1 -q_2)(q_1^2 -3q_1q_2 + q_2^2) = 0,$$ which is a (3,3)-net, has 9 lines, 10 triple points and 6 nodes. Such (3,3)-nets are discussed in \cite[Theorem 2.2 (3) and Example 2.3]{D2}.
\end{remark}

The next pencil is briefly discussed in \cite{Val}, see the Pappus arrangement completed in the last section. 

\begin{theorem}\label{thm2}
Consider the pencil of cubic curves given by $\PPP : u\CC_1 + v\CC_2$ in the complex projective plane $\PP^2_\C,$ where 
\begin{center}
$\CC_1 : q_1(x,y,z) = y(x-y-z)(2x+y-z) = 0$ and $\CC_2 : q_2(x,y,z) = xz(2x-5y+z) = 0.$
\end{center}

Then the following hold. 

\begin{enumerate}
\item The two curves $\CC_1$ and $\CC_2$ meet transversally in 9 points, hence the generic member of the pencil $\PPP$ is smooth.
\item The number $m$ of singular members in the pencil $\PPP$ is 6 and they correspond to the points $(u':v')\in \PP^1_\C$ that are listed as follows: for $(u' : v') = (1: 0), (u' : v') = (0 : 1)$ and $(u' : v') = (1 : 1)$ we get each time a triangle, and for $(u' : v') = (t_i : -1)$ where $i = 1,2,3$ and $t_i$'s are the roots of the equation $$125t^3 +399t^2 + 339t+1 = 0,$$ we get each time a nodal cubic.
\item The union of the singular members in the pencil $\PPP$ is a free curve $$\CC : f = q_1q_2(q_1 -q_2)(125q_1^3 -399 q_1^2 q_2 + 339q_1q_2^2 -1) = 0$$ of degree 18 and with exponents $d_1 = 4$ and $d_2 = 13,$ having 9 lines and three nodal cubics as irreducible components.
\item The monodromy operator (\ref{eq1}) associated to the curve $\CC$ has eigenvalues which are roots of unity of order 6. More precisely, the Alexander polynomial of the curve $\CC$ is given by 
$$\Delta_\CC(t)=(t-1)^7(t^6-1)^4.$$
\end{enumerate} 
\end{theorem}

\begin{remark}\label{rk2} The curve arrangement $\CC$ in Theorem \ref{thm2} contains the line arrangement $$\A : q_1q_2(q_1 -q_2) = 0,$$ which is a (3,3)-net, has 9 lines, 9 triple points and 9 nodes. Such (3,3)-nets are discussed in \cite[Theorem 2.2 (2) and Example 2.3]{D2}. 
\end{remark}

\begin{remark}\label{rk3}
All the irreducible components of the curve $\CC$ in Theorem \ref{thm1} are rational curves, and all the singularities of $\CC$ are ordinary multiple points, namely eight of order 2, one of order 3 and nine of order 5. In other words, our curve $\CC$ is very close to a line arrangement, both globally and locally. Note that the sum of the total Milnor numbers of the singular members in $\PPP$ is $$1 + 1 + 4 + 3 + 3 = 12 = 3(3-1)^2,$$ as predicted by the theory, see \cite{D1}, \cite{Val}. For the curve $\CC$ in Theorem \ref{thm2}, similar results apply, in particular in this case the sum of the total Milnor numbers of the singular members in $\PPP$ is $$1 + 1 + 1 + 3 + 3 + 3 = 3(3-1)^2.$$ All the singularities are weighted homogeneous for obvious reasons, fact which allows us to use the results from \cite{D1} and \cite{Val}.\\

The Alexander polynomials of the free curves $\CC$ in Theorem \ref{thm1} and Theorem \ref{thm2} can be computed using the algorithm described in \cite{D3}.
\end{remark}

\begin{remark}The characteristic polynomials $\Delta_\CC^q(t)$ of the operators $h^q, \,0\leq q \leq 2,$ are related by the following formula (see \cite[Proposition 4.1.21]{DS})
$$\Delta_\CC^0(t)\Delta^1_\CC(t)^{-1}\Delta^2_\CC(t)=(t^d-1)^{\chi(U)},$$
where $\chi(U)$ denotes the Euler characteristic of $U.$ Since the singularities are isolated, $\X(U)$ can be easily computed (see for instance \cite[Corollary 5.4.5]{DS}) and since the curve $\CC$ is reduced, $\Delta^0_\CC(t)=t-1.$ It follows that the operator $h^2$ is completely determined by $h^1$ and we can reduce to study the Alexander polynomial $\Delta_\CC(t):=\Delta_\CC^1(t).$
\end{remark}

 \section{Proof of Theorem \ref{thm1}}

\begin{enumerate}
\item Without loss of generality, one can assume that $x=1,$ and we have that $(1:y:z)$ is solution of 
the system
\begin{equation}\label{S0}
S:\left \{
\begin{array}{ccc}
    3yz + y^3 +z^3 & =0 \\
    3yz+ 1 + z^3 & =0 
\end{array}
\right.
\end{equation}
implies $y\in\{1,j,j^3\},\,j=\exp(2\sqrt{i}\pi/3),$ by subtracting the two equations. If $y=1$ then $z$ is solution of the equation 
$E_1: z^3+3z+1=0;$ if $y=j$ then $z$ is solution of the equation 
$E_2: z^3+3jz+1=0;$ and if $y=j^2$ then $z$ is solution of the equation 
$E_3: z^3+3j^2z+1=0.$ Denote by $\mathcal{S}_i$ the solution set of $E_i,$ with $i=1,2,3.$ It is clear that $\mathcal{S}_i \cap \mathcal{S}_j = \emptyset,$ for all $i \neq j$ and the system (\ref{S0}) has nine solutions:
$$\{(1:1:z),\,z\in \mathcal{S}_1\} \cup \{(1:j:z),\,z\in \mathcal{S}_2\} \cup \{(1:j^2:z),\,z\in \mathcal{S}_3\}.$$
 Hence $|\CC_1 \cap \CC_2| = 9$ and the intersection points of the curves $\CC_1$ and $\CC_2$ have intersection multiplicity 1, i.e. the two curves meet transversally and the base locus of the pencil $\PPP$ is smooth.

\item First recall that the member in the pencil $\mathcal{P}$ corresponding to the point $(u':v')\in\PP_{\C}^1$ is the curve 
$$u\mathcal{C}_1+v\mathcal{C}_2,$$
where $u=-v'$ and $v=u'.$ If $(u':v')=(1:0),$ the corresponding curve $\mathcal{C}_2$ is an irreductible nodal curve with node $(0:1:0).$ 
Similarily, for $(u':v')=(0:1)$ we get an irreductible nodal curve $\CC_1$ with node $(1:0:0).$ 
Suppose now that $u'$ and $v'$ are both non zero. Without loss of generality, one can assume $u'=1$ and we have that the corresponding member $u\CC_1+ \CC_2$ is singular if and only if the following system admits a solution:
\begin{equation}\label{S1}
\left \{
\begin{array}{ccc}
    3(u+1)xyz + uy^3 +(u+1)z^3 + x^3 & =0 & (i)\\
    (u+1)yz+ x^2 & =0 & (ii)\\
    (u+1)xz+uy^2 & =0 & (iii)\\
    (u+1)xy + (u+1)z^2& =0 &(iv)
\end{array}
\right.
\end{equation}

If $(u':v')=(1:1),$ then $-\CC_1+\CC_2: x^3-y^3$ is the union of three concurrent lines: 
\begin{center}
$x=y,x=jy$ and $x=j^2y.$
\end{center}
Assume now $(u':v')=(1:v'),\,v' \neq 1.$ The case $y=0$ is an absurd, since then $x=z=0$ with (ii) and (iv). If $y = 1,$ we have $x=0 \Rightarrow z=0$ with (ii) and $u=0$ with (i), while $z=0 \Rightarrow x=0$ with (ii) and $u=0$ with (i), which are both in contradiction with our hypothesis $v'\neq 0$. Finally, $(x:1:z),\,x,z \neq 0,$ is a solution of (\ref{S1}) if and only if 
\[
\left \{
\begin{array}{ccc}
    3(u+1)xz + u +(u+1)z^3 + x^3 & =0 & (i')\\
    (u+1)z+ x^2 & =0 & (ii')\\
    (u+1)xz+u & =0 & (iii')\\
     x + z^2& =0 &(iv')
\end{array}
\right.
\]
and we have with $(iv')$ and $(iii')$ that $z^3=\frac{u}{u+1},$  which implies with (ii') that $u^3+3u+1=0.$ Denote by $t_i,\,i=1,2,$ the solutions of the equation $$t^2+3t+1=0,$$ and by $\alpha_i$ a complex cubic root of $\frac{t_i}{t_i+1}.$ Then the singularities of the member $t_i\CC_1 + \CC_2$ associated to the point $(1:-t_i)$ are the following three points 
\begin{center}
$(-\alpha_i^2:1:\alpha_i),(-j^2\alpha_i^2:1:j \alpha_i)$ and $(-j\alpha_i^2:1:j^2\alpha_i),$
\end{center}
which are not colinear and hence form a triangle.
\item Since $t_1+t_2=-3$ and $t_1t_2=1,$ the union $\mathcal{C}$ of the singular members in the pencil $\PPP$ listed in 2. is the curve of degree $d=15$ defined by the homogeneous polynomial $$f = q_1q_2(q_1- q_2)(q_1^2-3q_1q_2 + q_2^2).$$ Since the singularities are all weighted homogeneous, we can apply the results from \cite{D1} (see Theorem 1.14) or \cite{Val} (see Theorem 2.7) to show that $\CC$ is free with exponents (4,10).
\item Recall that the $1-$eigenspace $H^1(F,\C)_1$ of the monodromy (\ref{eq1}) is a pure Hodge structure of type $(1,1),$ and that the sum of the non trivial eigenspaces $\displaystyle{\sum_{\lambda^d=1,\,\lambda \neq 1}H^1(F,\C)_\lambda }$ is pure of weight 1. For any $k\in [1,d],$ it is known that there exists a spectral sequence $E_*(f)_k$ whose first term is constructed from the Koszul complex in $\C[x,y,z]$ of the partial derivatives of $f,$ and whose limit $E_{\infty}(f)_k$ gives the action of the monodromy on the graded pieces $H^*(F,\C)_\lambda,\,\lambda=\exp(-2\sqrt{-1}\pi k/d),$ with respect to the pole order filtration $P,$ which contains the Hodge filtration $F$ and satisfies $P^2=0.$ In \cite{D3},  A. Dimca and G. Sticlaru showed that the computation of the second terms given in Equation (\ref{graduate}) is sufficient to detect all the non trivial eigenpaces of the operator (\ref{eq1}), see \cite[Theorem 1.2]{D3}. More in detail, for any $k\in [1,d]$, these second terms are given by

\begin{equation}\label{graduate}
E_2^{1,0}(f)_k= Gr_P^1 H^1(F,\C)_{\lambda}, 
\end{equation}
where $\lambda= \exp(-2\sqrt{-1}\pi k/d).$ Furthermore, when the curve is free, the authors also describe an algorithm which computes the dimensions $\dim E_2^{1,0}(f)_k.$ By using this algorithm and the computer algebra software Singular \cite{Sing}, we get that the only non zero dimensions of second terms of the form (\ref{graduate}) are listed as follows: $\dim E_2^{1,0}(f)_6=1,\,\dim E_2^{1,0}(f)_9=2,\,\dim E_2^{1,0}(f)_{12}=3,$ and $\dim E_2^{1,0}(f)_{15}=10.$  First, since $H^1(F,\C)_1=H^1(U,\C)$ and the curve $\mathcal{C}$ is the union of $11$ irreductible curves, it is well known that $\dim H^1(F,\C)_1=b_1(U)=10,$ see for instance \cite[Example 4.1]{D}. Using \cite[Theorem 1.2]{D3} we have that the multiplicities $m(\lambda)$ of the eigenvalues $\lambda \neq 1$ of the monodromy operator (\ref{eq1}) satisfy:
\begin{center}
$2 \leq m(\lambda)\leq 3$, for $k=6,9;$

$m(\lambda)= 3$, for $k=3,12;$ 

$m(\lambda)= 0$, for $k\notin \{3,6,9,12,15\}.$
\end{center}

The non trivial monodromies listed above are the roots of the unity of order 5. Since the monodromy operator $h^1$ is definied over $\Q,$ it is known that $\Delta_{\CC}(t) \in \Q[t]$ is a product of cyclotomic polynomials $\varphi_n$ with $n$ dividing $d=15$, and  it follows that 
$$\Delta_\CC(t)= (t-1)^7(t^5-1)^3.$$
\end{enumerate}
 
\section{Proof of Theorem \ref{thm2}}

\begin{enumerate}
\item Let us compute explicitely the intersection $\mathcal{C}_1\cap \mathcal{C}_2$ by solving the system of equations $q_1(x,y,z)=q_2(x,y,z)=0.$
If $x=0,$ then we get $(0:0:1),(0:1:1)$ and $(0:-1:1).$ Assume now $x=1.$ The first equation gives $y=0,y=1-z$ or $y=z-2.$
\begin{enumerate}
\item If $y=0,$ then $z=0$ or $-2$ and we get $(1:0:0)$ and $(1:0:-2);$
\item If $y=1-z,$ then $z=0$ or $\frac{1}{2}$ gives $(1:1:0)$ and $(2:1:1);$
\item If $y=z-2,$ then with $z=0$ or $z=3$ we find $(1:-2:0)$ and $(1:1:3).$
\end{enumerate}
Hence $|\CC_1 \cap \CC_2| = 9$ and $\CC_1$ and $\CC_2$ meet transversally.

\item Let us list all the singular members $u\mathcal{\CC}_1+v \mathcal{\CC}_2$ associated to points $(u':v')\in \PP_\C^1,$ where $u'=v$ and $v'=-u.$ 
Let us study separately the cases $(u':v')=(1:0)$ (a), $(u':v')=(0:1)$ (b), $(u':v')=(1:1)$ (c), and $(u':v')=(u':-1)$ with $u' \neq 0,-1$ (d).
\begin{enumerate}
\item[(a)] If $(u':v')=(1:0),$ then the singularities of the associated member $\mathcal{C}_2$ are the solutions of the sytem
\begin{equation}\label{S5}
\left \{
\begin{array}{ccc}
    xz(2x-5y+z) & =0 & (i)\\
   z(4x-5y+z) & =0 & (ii)\\
   -5xz & =0 & (iii)\\
   x(2x-5y+2z) & =0 & (iv) \\
\end{array}
\right.
\end{equation}
Without loss of generality, one can assume $y=1$ and we have $x=0$ or $z=0$ with $(iii).$ If $x=0,$ then $(ii)\Rightarrow z=0$ or $z=5$ that is $(0:1:0)$ and $(0:1:5)$ are solutions of (\ref{S5}). The case $z=0$ and $x=5/2$ implied by $(iv)$ gives the third solution $(5:2:0)$ and we get a triangle.

\item[(b)] Similarily, for $(u':v')=(0:1)$ the singularities of $\mathcal{C}_1$ are the solutions of 
\begin{equation}\label{S6}
\left \{
\begin{array}{ccc}
    y(x-y-z)(2x+y-z)& =0 & (i)\\
   y(4x-y-3z) & =0  & (ii)\\
   2x^2-3y^2-2xy-3xz+z^2& =0 & (iii)\\
   -y(3x-2z)& =0 & (iv)\\
\end{array}
\right.
\end{equation}

Without loss of generality, one can assume $x=1$ and we have with $(iv)$ that $y=0$ or $z=\frac{3}{2}.$ If $y=0,$ then $z$ is solution of $z^2-3z+2$ with $(iii)$ and $(1:0:1)$ and $(1:0:2)$ are solutions of (\ref{S6}). Assume  $z=\frac{3}{2}.$ Then $(i)$ and $(ii)$ both imply $y=0$ or $y=-\frac{1}{2},$ while equation $(iii)$ implies $y=-\frac{1}{2}$ or $-\frac{1}{6}.$ Hence $y=-\frac{1}{2}$ and $(2:-1:3)$ is the third solution of (\ref{S6}), i.e. we get once again a triangle.

\item[(c)]If $(u':v')=(1:1),$ then the singularities of the associated member $\mathcal{C}_1-\mathcal{C}_2$ are the solutions of the sytem

$$
\left \{
\begin{array}{cc}
    y(x-y-z)(2x+y-z)-xz(2x-5y+z) & =0  \\
   y(4x-y-3z)-z(4x-5y+z) & =0 \\
   2x^2-3y^2-2xy+2xz+z^2 & =0 \\
   -y(3x-2z)-x(2x-5y+2z) & =0 \\
\end{array}
\right.
$$
Without loss of generality, one can assume $y=1$ and the previous system becomes:
$$
\left \{
\begin{array}{ccc}
    (x-1-z)(2x+1-z)-xz(2x-5+z) & =0  & (i)\\
   4x-1-3z-z(4x-5+z) & =0 & (ii)\\
   2x^2-3-2x+2xz+z^2 & =0 & (iii)\\
   -3x+2z-x(2x-5+2z) & =0 & (iv)\\
\end{array}
\right.
$$
Hence $z(1-x)=x(x-1)$ with $(iv).$ If $x=1,$ then with equations $(ii)$ and $(iii)$ we have that $z^2+2z-3=0,$ which gives the solutions $(1:1:1)$ and $(1:1:-3).$ If $x\neq 1,$ then $(iv)\Rightarrow z=-x$ and we get the system
$$
\left \{
\begin{array}{cc}
    (2x-1)(3x+1)+x^2(x-5) & =0  \\
   3x^2+2x-1 & =0 \\
   x^2-2x-3 & =0 \\
   z & =-x \\
\end{array}
\right.
$$
with $(-1:1:1)$ as unique solution.
\item[(d)]Assume $(u':v')=(u':-1)$ with $u'\neq 0,-1.$ Then the associated curve $\mathcal{C}_1+v\mathcal{C}_2$ is singular whenever the following system has a solution:
$$
\left \{
\begin{array}{cc}
    y(x-y-z)(2x+y-z)+vxz(2x-5y+z) & =0  \\
   y(2x+y-z)+2y(x-y-z)+vz(2x-5y+z)+2vxz & =0 \\
   (x-2y-z)(2x+y-z)+y(x-y-z) - 5vxz & =0 \\
   -y(2x+y-z)-y(x-y-z)+vx(2x-5y+z) +vxz & =0 \\
\end{array}
\right.
$$
\end{enumerate}
Without loss of generality, one can assume that $y=1,$ that is:
$$
\left \{
\begin{array}{ccc}
    (x-1-z)(2x+1-z)+vxz(2x-5+z) & =0  & (i)\\
   (2x+1-z)+2(x-1-z)+vz(2x-5+z)+2vxz & =0 & (ii)\\
   (x-2-z)(2x+1-z)+(x-1-z) - 5vxz & =0 & (iii)\\
   -(2x+1-z)-(x-1-z)+vx(2x-5+z) +vxz & =0 & (iv)\\
\end{array}
\right.
$$
First, we have with $(iv)$ that 
\begin{equation}\label{eqz}
2z(1+vx)=x(-2vx+5v+3).
\end{equation}

If $1+vx=0,$ then $x=\frac{5v+3}{2v}$ with the previous equality. Since $x=-\frac{1}{v}$ we get $v=-1,$ which contradicts our assumption. Hence one can assume $1+vx\neq 0$ and by rewriting the equations of the previous system as polynomials in $(\C[x])[z]$ we get:
$$
\mathcal{S}:\left \{
\begin{array}{ccc}
    (1+vx)z^2 + x(-3-5v+2vx)z + 2x^2-x-1 & =0  & (i)\\
   vz^2 + (-3-5v+4vx)z +4x-1& =0 & (ii)\\
   z^2 + x(-3-5v)z +2x^2-2x-3 & =0 & (iii)\\
   \frac{1}{2(1+vx)}x(-2vx+5v+3)& =  z & (iv)\\
\end{array}
\right.
$$
Now, by injecting (\ref{eqz}) in $(i),$ the latter gives 
$$
\begin{array}{ccc}
    \frac{1}{1+vx} (2x^2-x-1) & =z^2  & (i')\\
\end{array}
$$
Let us now replace $z$ and $z^2$ by $(iv)$ and $(i')$ in $(ii)$ and $(iii)$ of $\mathcal{S},$ multiplicate $(ii)$ and $(iii)$ by $2(1+vx),$ and rewrite the obtained equations $(ii')$ and $(iii')$ as polynomials in $(\C[u])[x].$ Then the system $\mathcal{S}$ is equivalent to the following one:

$$
\mathcal{S'}:\left \{
\begin{array}{ccc}
    \frac{1}{1+vx} (2x^2-x-1) & =z^2  & (i')\\
   -8v^2x^3 + 30v(v+1)x^2 + (-25v^2-34v-1)x+ 2(-v-1)& =0 & (ii')\\
   10v(v+1)x^3 + (-25v^2-34v-1)x^2 + 6(-v-1)x -8& =0 & (iii')\\
   \frac{1}{2(1+vx)}x(-2vx+5v+3)& =  z & (iv)\\
\end{array}
\right.
$$
Remark that $(i')$ and $(iv)$ are equivalent whenever equations $(ii')$ and $(iii')$ are both satisfied. Indeed, from  $(i')$ and $(iv)$ we have that $$\frac{1}{4(1+vx)^2}x^2(-2vx+5v+3)^2= \frac{1}{1+vx}(2x^2-x-1),$$ which is equivalent to the following equation:
\begin{equation}\label{eq4}
4v^2x^4 -20v(1+v)x^3 + (25v^2+34v+1)x^2 + 4(1+v)x+4=0.
\end{equation}
Finally, by adding to equation $(iii')$ two times equation (\ref{eq4}), we get exactly equation $(ii'),$ up to multiplication by $(-x).$ It follows that the solutions of $\mathcal{S'}$ are given by the roots of the resultant of the two polynomials of $(ii')$ and $(iii')$. 
One can compute that this resultant is 
$$110592v^3(v+1)^3(125v^3+399v^2+339v+1)$$
and if $t_i,\,i=1,2,3,$ are the roots of the equation
$$125t^3+399t^2+339t+1,$$
then the curve $\mathcal{C}_1+t_i\mathcal{C}_2$ is singular with singular points $(x:1:z),$ where $x$ is solution of 
$$
\left \{
\begin{array}{ccc}
   -8t_i^2x^3 + 30t_i(t_i+1)x^2 + (-25t_i^2-34t_i-1)x+ 2(-t_i-1)& =0 & (ii')\\
   10t_i(t_i+1)x^3 + (-25t_i^2-34t_i-1)x^2 + 6(-t_i-1)x -8& =0 & (iii')\\
\end{array}
\right.
$$
 and $z=\frac{1}{2(1+t_ix)}x(-2t_ix+5t_i+3).$
By performing $4(ii') + (t_i-1)(iii')$ we have that $x$ is solution of 
$$
  2t_i(5t_i-1)(t_i-5)x^2-(t_i-1)(25t_i^2-154t_i+1)x -2(47t_i^2-62t_i-1)=0.
$$
Then it is possible to guess the number of singularities of the three curves $\mathcal{C}_1+t_i\mathcal{C}_2,\,i=1,2,3,$ by using a cardinality argument and the total Milnor number of the curve $\mathcal{C}.$ Indeed, since the pencil is generic from 1., we have with \cite[Proposition 5.1]{D1} that the sum of the Milnor numbers of all the singularities of the degree $3$ members $u\CC_1+v\CC_2$ listed before is equal to $3(3-1)^2=12.$ It follows that the not yet known singularities (at least one for each curve $\CC_1+t_j\CC_2$) contribute to $12-(3\times 1 + 3\times 1+ 3\times 1)=3.$ Hence each member $\CC_1+t_i\CC_2$ has exactly one singularity and is an irreductible nodal curve. 
\item Since $t_1+t_2+t_3=\frac{-399}{125},\,t_1t_2+t_1t_3+t_2t_3=\frac{339}{125}$ and $t_1t_2t_3=\frac{-1}{124}$ with Viete formula, the union $\CC$ of the singular members of the pencil $\mathcal{P}$ listed above is defined by the homogeneous degree 18 polynomial
$$f= q_1q_2(q_1-q_2)(125q_1^3-399q_1^2q_2 + 339q_1q_2^2 -1).$$
Since the singularities are all weighted homogeneous, we can deduce from the results in \cite{D1} and \cite{Val} that $\CC$ is free with exponents (4,13).

\item
By applying the algorithm described in \cite{D3} we get this time that the only non zero dimensions of the second terms of the form (\ref{graduate}) are listed as follows: $\dim E_2^{1,0}(f)_6=1,\, \dim E_2^{1,0}(f)_9=2,\, \dim E_2^{1,0}(f)_{12}=3,\, \dim E_2^{1,0}(f)_{15}=4,$ and $\dim E_2^{1,0}(f)_{18}=11.$
Using \cite[Theorem 1.2]{D3} we have that the multiplicities $m(\lambda)$ of the eigenvalues $\lambda \neq 1$ of the monodromy operator (\ref{eq1}) satisfy:
\begin{center}
$m(\lambda)=4,$ for $k=3,15$ ($\lambda$ of order 6);\\
$3 \leq m(\lambda)\leq 4$, for $k=6,12$ ($\lambda$ of order 3);\\
$2 \leq m(\lambda)\leq 4$, for $k=9$ ($\lambda$ of order 2);\\
$m(\lambda)= 0$, for $k\notin \{3,6,9,12,15,18\}.$
\end{center}
In particular, $H^1(F,\C)_\lambda \neq 0 \Rightarrow \dim H^1(F,\C)_\lambda \leq 4,$ if $\lambda \neq 1.$

Let $S:= \PP^1_\C \backslash \{(u':v') \in \PP^1_\C \,\,|\,\,(u':v') \text{\,is\,a\,singular\,member\,in\,}\PPP\}$ be the complement of the six points described in 2. Then by considering the surjective morphism $r: U \rightarrow S,\,r(x:y:z)=(q_1(x,y,z):q_2(x,y,z))$, and applying
\cite{Ar}, \cite[Theorem 3.1]{DP} or \cite[Corollary 4.7]{D0} we get that 
$$\dim H^1(F)_{\lambda^k} \geq -\chi(S)=4,$$
for any $\lambda^k \neq 1.$
Hence the Alexander polynomial has de following expression:
$$\Delta_\CC(t)=(t-1)^7(t^6-1)^4.$$
\end{enumerate}

\medskip

\noindent
{\bf Acknowledgements.}  
The author thanks Alexandru Dimca for his 
suggestions and motivating discussions which led her to obtain the results in this paper. She is supported by the BREMEN TRAC Postdoctoral Fellowships for Foreign Researchers.


\begin{thebibliography}{99}
\bibitem{Ar} D. Arapura: Geometry of cohomology support loci for local systems. I, J. Algebraic Geom. {6} (1997), 563–597.
\bibitem{Sing}
W. Decker, G.-M. Greuel, G.Pfister, H. Sch{\"o}nemann. \newblock {\sc Singular} {4-0-2} --- {A} computer algebra system for polynomial computations.
\newblock {http://www.singular.uni-kl.de} (2015).
\bibitem{D0} A. Dimca: Characteristic varieties and constructible sheaves, Rend. Lincei Mat. Appl. 18 (2007), 365- 389.
\bibitem{D1}A. Dimca: Curve arrangements, pencils, and Jacobian syzygies, arXiv: 1601:00607.
\bibitem{D} A. Dimca: On the Milnor fibrations of weighted homogeneous polynomials, Compositio Math. 76 (1990), 19-47.
\bibitem{DS} A. Dimca: Singularities and topology of hypersurfaces. Universitext. Springer-Verlag, New York, 1992.
\bibitem{D2}A. Dimca, D. Ibadula, D. A. M\v{a}cinic: Pencil type line arrangements of low degree: classiﬁcation and monodromy, Ann. Sc. Norm. Super. Pisa XV(2016),249-267.
\bibitem{D3}A. Dimca, G. Sticlaru: A computational approach to Milnor ﬁber cohomology, arXiv:1602.03496.

\bibitem{DP} A. Dimca, S. Papadima: Finite Galois covers, cohomology jump loci, formality properties, and multinets,  Ann. Scuola Norm. Sup. Pisa Cl. Sci (5), Vol. X (2011), 253-268. 

\bibitem{Falk}M. Falk, S. Yuzvinsky: Multinets, resonance varieties, and pencils of plane curves, Compositio Math. 143 (2007), 1069–1088.

\bibitem{PS}S. Papadima, A. I. Suciu: The Milnor fibration of a hyperplane arrangement: from modular resonance to algebraic monodromy, arXiv:1401.0868.

\bibitem{ST}H. K. Schenck, S.O. Toh\v{a}neanu: Freeness of conic-line arrangements in $\PP^2$, Comment. Math. Helv. 84(2009), 235-258.
\bibitem{Val}J. Vall\`{e}s: Free divisors in a pencil of curves, Journal of Singularities 11 (2015), 190-197, and Erratum: Journal of Singularities 14 (2016), 1-2.
\bibitem{Yuz}S. Yuzvinsky: A new bound on the number of special ﬁbers in a pencil of curves. Proc. Amer. Math. Soc. 137 (2009), 1641–1648.

\end{thebibliography}
\end{document}